\documentclass{amsart}

\usepackage{amssymb}
\usepackage{amsmath,amscd}
\usepackage{amsthm}
\usepackage{graphicx,psfrag}

\theoremstyle{plain}
\newtheorem{theorem}{Theorem}[section]
\newtheorem{thm}{Theorem}

\newtheorem{lemma}[theorem]{Lemma}

\theoremstyle{definition}
\newtheorem{definition}[theorem]{Definition}

\theoremstyle{remark}
\newtheorem{remark}[theorem]{Remark}
\newtheorem*{notation}{Notation}
\newtheorem*{acknowledgments}{Acknowledgments}
\newtheorem{property}{Property}
\newtheorem{step}{Step}
\newtheorem{observation}{Observation}
\newtheorem{case}{Case}
\newtheorem{S4case}{Case}
\newtheorem{subcase}{Subcase}[case]
\newtheorem{S4subcase}{Subcase}[S4case]

\numberwithin{figure}{section}

\begin{document}

\title{An algorithm to find vertical tori in small Seifert fiber spaces}

\author{Tao Li}
\thanks{Partially supported by an NSF grant} 

\address{Department of Mathematics \\
 Oklahoma State University \\
 Stillwater, OK 74078}
\email{tli@math.okstate.edu}

\begin{abstract}
We give an algorithm to find vertical essential tori in small Seifert fiber spaces with infinite fundamental groups.  This implies that there are algorithms to decide whether a 3-manifold is a Seifert fiber space.  
\end{abstract}
\maketitle

\begin{psfrags}
\psfrag{(a)}{(a)}
\psfrag{(b)}{(b)}
\psfrag{(c)}{(c)}
\psfrag{(d)}{(d)}
\psfrag{(e)}{(e)}
\psfrag{(f)}{(f)}
\psfrag{A}{$A$}\psfrag{A'}{$A'$}
\psfrag{B}{$B$}\psfrag{B'}{$B'$}
\psfrag{E}{$E$}
\psfrag{O}{$O$}
\psfrag{P}{$P$}
\psfrag{X}{$X$}
\psfrag{Y}{$Y$}
\psfrag{Z}{$Z$}
\psfrag{tau}{$\tau$}
\psfrag{tauA}{$\tau_A$}
\psfrag{eta}{$\eta$}
\psfrag{etaA}{$\eta_A$}
\psfrag{etaB}{$\eta_B$}
\psfrag{delta}{$\delta$}

\section{Introduction}
A fundamental problem in 3-manifold topology is to recognize a 3-manifold from a given combinatorial structure, e.g., a triangulation or a handle decomposition. Much progress has been made in the past few years.  In the 1980s, Jaco and Oertel \cite{JO} showed that there is an algorithm to decide whether an irreducible 3-manifold is Haken.  Later, Rubinstein and Thompson gave algorithms to recognize 3-spheres \cite{R, Th}, indicating that there are also algorithms to decide whether a 3-manifold is reducible.  Recently, it was shown in \cite{AL} that there are algorithms to determine whether a $3$-manifold contains essential laminations or taut foliations.  Seifert fiber spaces are a major class of 3-manifolds and have important roles in Thurston's geometrization conjecture.  In this paper, we will give an algorithm to determine whether a 3-manifold is a Seifert fiber space.  

\begin{thm}\label{T1}
There is an algorithm to recognize Seifert fiber spaces.
\end{thm}

Let $M$ be a 3-manifold.  Since there is an algorithm to decide whether a manifold is Haken \cite{JO}, and since there is an algorithm to recognize Haken Seifert fiber spaces \cite{JT}, we only need to consider small Seifert fiber spaces.  Using strongly irreducible Heegaard splitting and almost normal surfaces, Rubinstein gave an algorithm to recognize lens spaces \cite{R}.  Moreover, Rubinstein's techniques also give an algorithm to decide whether a 3-manifold is a small Seifert fiber space with finite fundamental group.  So, the key part of Theorem~\ref{T1} is to find an algorithm to recognize small Seifert fiber spaces with infinite fundamental groups.  Peter Scott has shown that each small Seifert fiber space with infinite $\pi_1$ contains a vertical essential torus with only one or two double curves and having the 4-plane property \cite{S}.  The main goal of this paper is to give an algorithm to find such a vertical torus.  Once we find such a torus, we immediately see the Seifert fiber structure by checking the complement of this immersed vertical torus.

In this paper, we first show that the essential tori considered in \cite{S} have the 7-color property, i.e., one can use $7$ different colors to color their preimages in the universal cover so that any two planes in the same color do not intersect each other.  Using this result, we are able to construct finitely many immersed branched surfaces, one of which fully carries such an immersed torus.  Then, we analyze these immersed branched surfaces and get an upper bound on the weight of some essential torus.  One can easily enumerate all the immersed normal tori with weight bounded by this number, and by a theorem in \cite{LN}, one of these tori is an essential torus with the $4$-plane and $1$-line properties.  For any essential torus with the $4$-plane and $1$-line properties, one can perform some simple homotopies to eliminate all triple points and get a vertical torus.  So, we can algorithmically perform such homotopies on each of these tori.  If we cannot eliminate the triple points by such simple homotopies for any torus, we can conclude that the 3-manifold is not a Seifert fiber space (with infinite $\pi_1$); if we get a torus without triple points, we only need to check whether its complement is a union of solid tori.

A different algorithm has also been found by Rannard and Rubinstein \cite{RR}.  I thank Dave Letscher for providing me a copy of \cite{RR}.  The method used in \cite{RR} is an analysis of a sequence of involutions on genus 2 Heegaard splittings.  Ian Agol has also independently obtained an algorithm by studying the representation of the fundamental groups and using computational algebraic geometry.  Although none of the three algorithms are efficient, the algorithm in this paper is much easier to implement.

\begin{acknowledgments}
Part of this paper was written while I was visiting the American Institute of Mathematics, and I would like to thank AIM for its hospitality.  I would also like to thank Hyam Rubinstein for pointing out a gap in an earlier version, and I thank Max Neumann-Coto for the discussions of triangle groups.
\end{acknowledgments}

\section{The n-color property}\label{S2}

\begin{definition}
Let $\Gamma$ be a collection of objects (usually lines or surfaces) in a manifold.  We say that $\Gamma$ satisfies the \emph{$n$-color property} if we can color all objects in $\Gamma$ using $n$ different colors such that any two objects in the same color do not intersect each other, i.e. $\Gamma$ can be divided into $n$ subsets $\Gamma_1,\dots,\Gamma_n$ such that for any $P, Q\in\Gamma_i$ ($i=1,\dots ,n$), $P\cap Q=\emptyset$.
\end{definition}

Let $M$ be a closed 3-manifold.  We consider $\pi_1$-injective surfaces in $M$.  By \cite{SU, SY}, there exists a least area map in the homotopy class of each $\pi_1$-injective surface.  Moreover, the preimage of any least area surface in the universal cover consists of embedded planes \cite{FHS}.  The same is also true for least weight immersed normal surfaces \cite{JR}.

\begin{definition}
We say that the immersed surface $f:S\to M$ (or simply $S$) has the \emph{$n$-color property} if $\pi^{-1}(f(S))$ is a set of planes in $\widetilde{M}$ satisfying the $n$-color property, where $\widetilde{M}$ is the universal cover of $M$ and $\pi:\widetilde{M}\to M$ is the covering map.  An immersed surface $f:S\to M$ has the \emph{$n$-plane property} \cite{HS} if $\pi^{-1}(f(S))$ is a set of planes and there is a disjoint pair among any collection of $n$ planes.
\end{definition}

\begin{remark}\label{R:color}
\begin{enumerate}
\item If $M$ has a Haken $n$-fold cover, then $M$ contains a surface with the $n$-color property.

\item If a surface has the $n$-color property, then it has the $(n+1)$-plane property.  
\item Rubinstein informed me that he and Sageev have shown that every closed $\pi_1$-injective surface in a hyperbolic 3-manifold satisfies the $k$-color property for some $k$.  They proved this for $k$-plane property in \cite{RS} earlier.  
\end{enumerate}
\end{remark}

\begin{notation}
Throughout this paper, we denote the number of components of $X$ by $|X|$, denote the interior of $X$ by $int(X)$.  We denote a group generated by $h_1,\dots,h_n$ by $<h_1,\dots,h_n>$.
\end{notation}

Let $M$ be a small Seifert fiber space with infinite $\pi_1$.  Since we will refer to \cite{S} many times, we keep the same notation.  The orbifold of $M$ is a 2-sphere $Q$ with 3 cone points $\bar{X}$, $\bar{Y}$ and $\bar{Z}$, with cone angles $2\pi/p$, $2\pi/q$, and $2\pi/r$ respectively, where $p,q,r$ are integers.  The fundamental group of $M$ has infinite cyclic center with quotient a hyperbolic or Euclidean triangle group $\Delta (p,q,r)$.  Let $\widetilde{Q}$ be the universal covering orbifold of $Q$, then $\widetilde{Q}$ has either Euclidean or hyperbolic structure depending on whether $1/p+1/q+1/r$ is equal to $1$ or less than $1$.  One can connect $\bar{X}$, $\bar{Y}$ and $\bar{Z}$ by 3 geodesic arcs in the induced Euclidean or hyperbolic metric.  The preimage of the 3 geodesic arcs gives a tessellation of the plane $\widetilde{Q}$ by geodesic triangles with angles $\pi/p$, $\pi/q$ and $\pi/r$.  Let $XYZ$ be a geodesic triangle in this tessellation of $\widetilde{Q}$, where $X$, $Y$ and $Z$ are in the preimage of $\bar{X}$, $\bar{Y}$ and $\bar{Z}$ respectively.  Then, the triangle group $\Delta (p,q,r)$ is generated by $x$, $y$ and $z$ which are clockwise rotations about $X$, $Y$, and $Z$ respectively through angles of $2\pi/p$, $2\pi/q$ and $2\pi/r$ ($xyz=1$).  

Let $\psi:M\to Q$ be the Seifert fibration.  We say an immersed torus $T$ is \emph{vertical} if $\psi^{-1}(\psi(T))=T$, i.e., $T$ consists of circle fibers of $M$.  Every loop $l$ in the orbifold $Q$ determines a vertical torus $\psi^{-1}(l)$, and $\psi^{-1}(l)$ is $\pi_1$-injective if and only if $l$ represents an element with infinite order in $\Delta (p,q,r)$.  In \cite{S}, Scott analyzes in depth the loops in $Q$ with one or two double points that represent $xy^{-1}$ or $xy^{-2}$, as shown in Figure~\ref{F31}.  He has shown:

\begin{lemma}[Scott \cite{S}]\label{Scott}
Suppose $\Delta(p,q,r)$ is a hyperbolic triangle group and $x$, $y$, $z$ are the rotations as above.  
\begin{enumerate}
\item If $p,q,r\ge 3$, $xy^{-1}$ has infinite order in $\Delta (p,q,r)$.
\item For $\Delta(p,q,2)$ ($p,q\ge 4$), $xy^{-1}$ has infinite order.
\item For $\Delta (3, q, 2)$, $q\ge 7$, $xy^{-2}$ has infinite order.
\end{enumerate}
\end{lemma}

In fact, Scott's lemma says the following. 

\begin{lemma}[Scott \cite{S}]\label{Scott2}
Let $x$ and $y$ be two elements in the group of orientation preserving isometries of $\mathbb{H}^2$ such that: 
\begin{enumerate}
\item $x$, $y$ and $xy$ are all of finite order, 
\item $x$ and $y$ generate a discrete infinite  group. 
\end{enumerate} 
Then, there is an element of infinite order that can be expressed by $uv^{\pm1}$ or $uv^{\pm 2}$, where $u$ and $v$ are among $x$, $y$, $xy$.
\end{lemma}
\begin{proof}

Since $x$, $y$ and $xy$ are all of finite order, $x$ and $y$ generate a hyperbolic triangle group, and the lemma follows from Lemma~\ref{Scott}.
\end{proof}

The vertical tori corresponding to the loops that represent $xy^{-1}$ or $xy^{-2}$ in Lemma~\ref{Scott} have one or two double curves.  Scott has shown that the preimage of these loops (representing $xy^{-1}$ or $xy^{-2}$ as in Lemma~\ref{Scott} and Figure~\ref{F31}) in the universal covering orbifold $\widetilde{Q}$ is a union of lines such that the intersection of each pair of lines consists of at most one point. The intersection patterns of the planes in the preimage of corresponding vertical tori (in the universal cover of $M$) is the same as the intersection patterns of the lines in $\widetilde{Q}$.  In this section, we will show that the preimages of these loops in $\widetilde{Q}$ have the 7-color property, and hence the corresponding vertical tori have the 7-color property.

\begin{figure}
\begin{center}
\psfrag{Xb}{$\bar{X}$}
\psfrag{Yb}{$\bar{Y}$}
\psfrag{Zb}{$\bar{Z}$}
\psfrag{xy1}{$xy^{-1}$}
\psfrag{xy2}{$xy^{-2}$}
\includegraphics[width=4in]{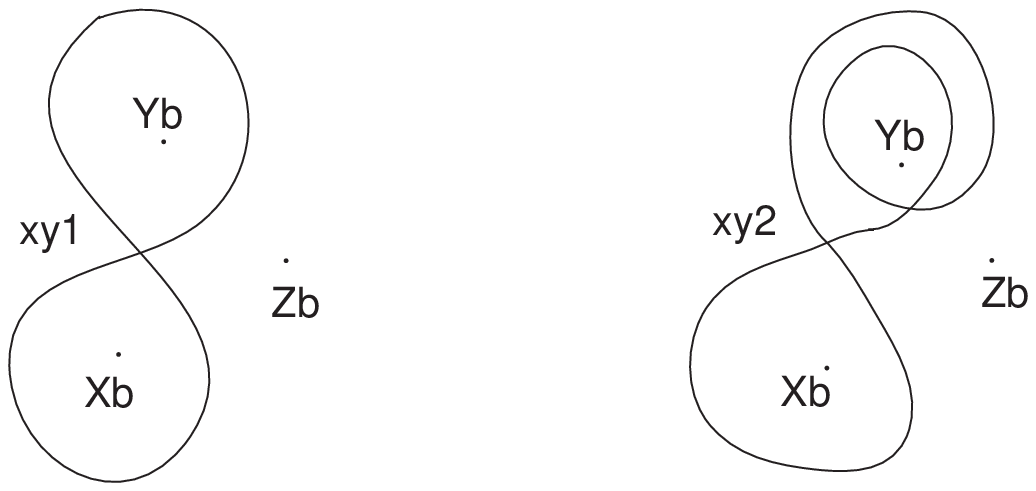}
\caption{}\label{F31}
\end{center}
\end{figure}

\begin{lemma}\label{L32}
If $M$ is a small Seifert fiber space of Euclidean type, then $M$ contains a vertical torus with the 3-color property.
\end{lemma}

\begin{proof}
Since $M$ is a small Seifert fiber space of Euclidean type, the triangle group of corresponding orbifold is either $\Delta(2,4,4)$, or $\Delta(2,3,6)$ or $\Delta(3,3,3)$.  If $M$ is of type $\Delta(2,4,4)$ (resp. $\Delta(2,3,6)$), then $M$ is double (resp. triple) covered by a Haken Seifert fiber space whose orbifold is a 2-sphere with 4 cone points.  If $M$ is of type $\Delta(3,3,3)$, then $M$ is a triple covered by a Haken Seifert fiber space whose orbifold is a torus without cone point.  Thus, the lemma follows from Remark~\ref{R:color} (1).
\end{proof}

\begin{theorem}\label{T31}
Let $M$ be a small Seifert fiber space with infinite $\pi_1$, then $M$ contains a vertical $\pi_1$-injective torus with the 7-color property.
\end{theorem}

\begin{proof}
By Lemma~\ref{L32}, we only need to consider the case that $M$ is a small Seifert fiber space of hyperbolic type.  Throughout the proof, we use $C$ to denote the loop representing $xy^{-1}$ or $xy^{-2}$ as in Lemma~\ref{Scott} and Figure~\ref{F31}, and use $\widetilde{C}$ to denote the preimage of $C$ in $\widetilde{Q}$.  Since $xy^{-1}$ (or $xy^{-2}$) has infinite order and $\widetilde{Q}$ is a hyperbolic plane, $xy^{-1}$ (or $xy^{-2}$) has an axis $l$ in $\widetilde{Q}$.  If $l$ does not pass through any translates of $X$, $Y$, or $Z$, then we can choose $C$ to be $\pi (l)$, where $\pi :\widetilde{Q}\to Q$ is the orbifold covering.  If $l$ passes through some translates of $X$, $Y$, or $Z$, let $m$ be a line of points that have a fixed distance $\epsilon$ from $l$ and on one side of $l$ for some small $\epsilon$, then $C$ can be chosen as $\pi (m)$.  In fact, by section 1 of \cite{S} (see Figures 5--9 in \cite{S}), $l$ passes through some translates of $X$, $Y$ or $Z$ if and only if $M$ is of type $\Delta (3,q,2)$ or $\Delta(4,q,2)$.   Let $\mathcal{L}$ be the collection of all the translates of $l$ in the hyperbolic plane $\widetilde{Q}$.  So, $\mathcal{L}$ has the 7-color property if and only if $\widetilde{C}$ has the 7-color property.  

We denote by $\mathcal{L}_0$ the union of intersection points of the geodesics in $\mathcal{L}$.  The geodesics in $\mathcal{L}$ give rise to a tiling of the hyperbolic plane $\widetilde{Q}$ with each tile the closure of a component of $\widetilde{Q}-\mathcal{L}$.  By section 1 of \cite{S}, there are no 3 geodesics in $\mathcal{L}$ passing through the same point in $\mathcal{L}_0$ in any case. So, $\mathcal{L}_0$ is the set of vertices in this tiling, and there are exact 4 tiles incident to each vertex.  If the geodesics in $\mathcal{L}$ do not pass through translates of $X$, $Y$ or $Z$, then its projection to the orbifold $Q$ is exactly as in Figure~\ref{F31}, and one can easily draw the tiling according to Figure~\ref{F31} and the indices of the cone points.  

We call an embedded disk $P$ in $\widetilde{Q}$ a \emph{polygon in this tiling} if $P$ is a polygon with vertices in $\mathcal{L}_0$ and each edge a geodesic arc (in a geodesic of $\mathcal{L}$).  Note that the interior of an edge of $P$ may contain other points in $\mathcal{L}_0$ and we may assume the inner angle of each vertex of $P$ is not $\pi$.  Thus, a polygon $P$ is convex if and only if the inner angle of every vertex of $P$ is less than $\pi$, and in particular, each tile in this tiling is a convex polygon.

This tilling can be considered a cellulation with $\mathcal{L}_0$ the 0-cells, $\mathcal{L}-\mathcal{L}_0$ the 1-cells, and $\widetilde{Q}-\mathcal{L}$ the 2-cells.  Let $P$ be a convex polygon in this tiling and $P''$ be the union of $P$ and all the tiles incident to $\partial P$.  Although we do not know whether $P''$ is a disk yet, we can construct a disk $P'$ by gluing corresponding tiles along $\partial P$, and construct a cellular map $f:P'\to\widetilde{Q}$ such that $f(P')=P''$ and $f$ is a local embedding (in particular, $f$ maps tiles to tiles).  So, we can consider $P'$ as a polygon ($P\subset int(P')$) with induced hyperbolic metric and induced tiling from $\widetilde{Q}$. The intersection of each tile in $P'-int(P)$ with $P$ is either a vertex or a single edge of this tile.  Next, we analyze the polygon $P'$.  

Let $\alpha_1,\dots,\alpha_k$ be those edges (of tiles in $P'-int(P)$) with exactly one endpoint in $\partial P$ (clockwise around $\partial P$).  Let $x_i$ and $y_i$ be the two endpoints of $\alpha_i$ with $x_i\in\partial P$.  By our construction $y_i\in\partial P'$.  Moreover, since our tiling comes from the intersection of geodesics and since $P$ is convex, $y_i$ lies in the interior of an edge of the polygon $P'$ unless $y_i=y_{i\pm 1}$, in which case $\alpha_i$, $\alpha_{i\pm 1}$ and the arc in $\partial P$ connecting $x_i$ and $x_{i\pm 1}$ bound a triangle tile in $P'-int(P)$, as shown in Figure~\ref{tiling}.  Thus, if there is no triangle tile in this tiling, the inner angle of each vertex of $P'$ is less than $\pi$ and hence $P'$ must also be a convex polygon.  This implies that $f:P'\to\widetilde{Q}$ is an embedding, i.e. $P''$ is also a convex polygon in the tiling.  Furthermore, suppose there is no triangle tile in this tiling and there is an edge $\gamma$ of $P'$ whose interior contains 2 points of $\mathcal{L}_0$.  Then, we can suppose $y_i\in int(\gamma)$ and $y_{i+1}\in int(\gamma)$ for some $i$.  Since $P$ is convex, $\alpha_i$, $\alpha_{i+1}$ must belong to a quadrilateral tile in $P'-int(P)$, and the other two edges of this quadrilateral tile lie in two edges of $P$ and $P'$ respectively.  Hence, if the interior of an edge of $P'$ contains 3 points of $\mathcal{L}_0$, then there must be two quadrilateral tiles sharing an edge.  Thus, we obtain the following observation.

\begin{observation}\label{O1}
If there is no triangle tile in this tiling by $\mathcal{L}$, then $P'$ is a convex polygon and $f:P'\to\widetilde{Q}$ is an embedding.  If, in addition, there are no two quadrilateral tiles sharing an edge in this tiling, then the interior of any edge of $P'$ contains at most two points of $\mathcal{L}_0$.  
\end{observation}

\begin{figure}
\begin{center}
\psfrag{P}{$P$}
\psfrag{P'}{$P'$}
\psfrag{e}{$e$}
\psfrag{xi}{$x_i$}
\psfrag{xi1}{$x_{i+1}$}
\psfrag{ali}{$\alpha_i$}
\psfrag{ali1}{$\alpha_{i+1}$}
\psfrag{vi}{$V_i$}
\psfrag{vi'}{$V_i'$}
\psfrag{vi''}{$V_i''$}
\psfrag{vi1}{$V_{i+1}$}
\psfrag{ti}{$t_i$}
\psfrag{ti1}{$t_{i+1}$}
\psfrag{ei'}{$e_i'$}
\psfrag{ei''}{$e_i''$}
\psfrag{si}{$s_i$}
\includegraphics[width=3in]{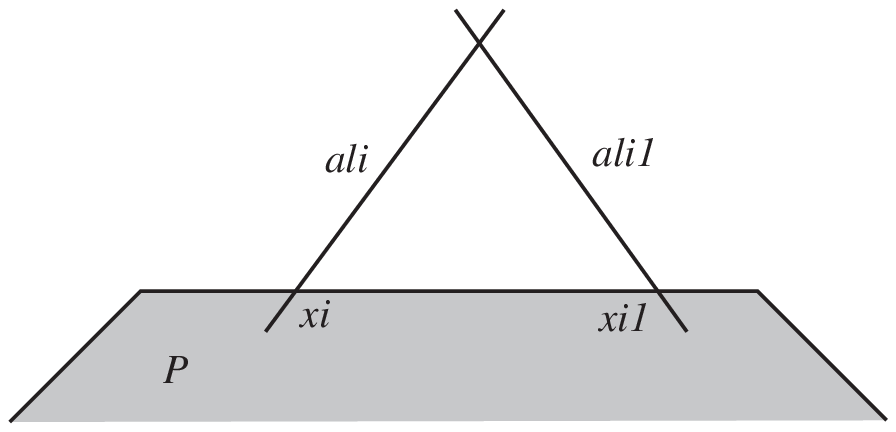}
\caption{} \label{tiling}
\end{center}
\end{figure}

\begin{figure}
\begin{center}
\psfrag{P}{$P$}
\psfrag{P'}{$P'$}
\psfrag{e}{$e$}
\psfrag{xi}{$x_i$}
\psfrag{xi1}{$x_{i+1}$}
\psfrag{ali}{$\alpha_i$}
\psfrag{ali1}{$\alpha_{i+1}$}
\psfrag{vi}{$V_i$}
\psfrag{vi'}{$V_i'$}
\psfrag{vi''}{$V_i''$}
\psfrag{vi1}{$V_{i+1}$}
\psfrag{ti}{$t_i$}
\psfrag{ti1}{$t_{i+1}$}
\psfrag{ei'}{$e_i'$}
\psfrag{ei''}{$e_i''$}
\psfrag{ei1'}{$e_{i+1}'$}
\psfrag{si}{$s_i$}
\psfrag{sigma}{$\sigma$}
\includegraphics[width=4in]{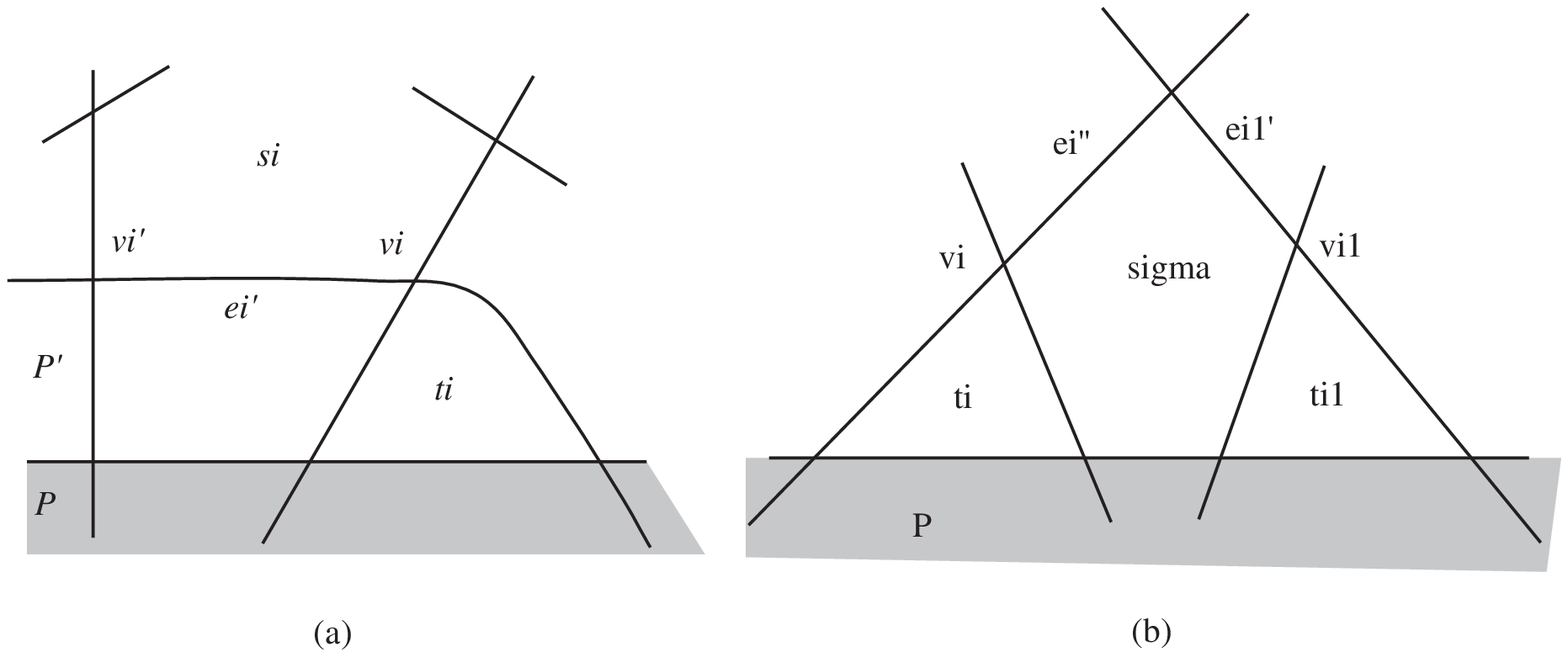}
\caption{} \label{triangle}
\end{center}
\end{figure}

\begin{case}\label{case1}
$M$ is not of type $\Delta(4,q,2)$ ($q>4$) or $\Delta(3,q,2)$ ($q>6$).  
\end{case}

In this case, the geodesics in $\mathcal{L}$ do not pass through translates of $X$, $Y$ or $Z$, and the loop $\pi(\mathcal{L})$ in the orbifold $Q$ is a figure eight representing $xy^{-1}$ as shown in Figure~\ref{F31}~(a).  Hence, the tiling from $\mathcal{L}$ consists of $p$-gons, $q$-gons and $2r$-gons, where $M$ is of type $\Delta(p,q,r)$.  The $p$-gons and $q$-gons only share edges with the $2r$-gons. 

Next two observations are easy to see from the indices of the cone points. 

\begin{observation}\label{O2}
If $M$ is not of type $\Delta(4,q,2)$ ($q>4$) or $\Delta(3,q,2)$ ($q>6$), there are no two quadrilateral tiles (in the tiling from $\mathcal{L}$) sharing an edge.
\end{observation}

\begin{observation}\label{O3} Suppose $M$ is not of type $\Delta(4,q,2)$ ($q>4$) or $\Delta(3,q,2)$ ($q>6$).  If an $n$-gon tile shares an edge with a triangular tile, then $2r=n\ge 6$.
\end{observation}

As before, if there is no triangular tile, then $f:P'\to\widetilde{Q}$ is an embedding and $f(P')$ is convex.  Now, we consider the case that there are triangle tiles in this tiling.  By the arguments on $P'$ before, $P'$ is not convex if and only if there is a triangle tile in $P'-int(P)$ with one edge in $\partial P$, as shown in Figure~\ref{tiling}.  Let $t_1,\dots,t_m$ be the collection of such triangle tiles (clockwise around $\partial P$), and let $V_i$ be the vertex of $t_i$ lying in $\partial P'$ for each $i$.  So, the inner angle (of $P'$) at $V_i$ is greater than $\pi$, i.e., there are 3 tiles in $P'-int(P)$ sharing the vertex $V_i$ for each $i$.  We use $s_i$ to denote the fourth tile incident to $f(V_i)$.  Similar to the construction of $P'$ from $P$, we can first glue a copy of $s_1$ to $P'$ to get a polygon $P_1'$ and extend the map $f:P'\to\widetilde{Q}$ to $f:P_1'\to\widetilde{Q}$ such that $f$ is a cellular map and a local embedding.  Suppose we have inductively constructed such a polygon $P_i'$ and extended the map $f$ to $f:P_i'\to\widetilde{Q}$.  If $V_{i+1}$ lies in the interior of $P_i'$ (we will show later that this is impossible), we let $P_{i+1}'=P_i'$ and $s_{i+1}=s_i$.  If $V_{i+1}\in\partial P_i'$, we can glue a copy of $s_{i+1}$ to $P_i'$ (according to the local picture of $f(P_i')$ at $f(V_{i+1})$) to get a polygon $P_{i+1}'$ and extend the map $f$ to a cellular map $f:P_{i+1}'\to\widetilde{Q}$ such that $f$ is a local embedding, in particular, two adjacent edges (of tiles) in $\partial P'_{i+1}$ are not mapped to the same 1-cell in $\mathcal{L}$.  In fact, we can construct $P_{i+1}'$ as follows. We first identify an edge (of a tile) in $\partial P_i'$ containing $V_{i+1}$ to the corresponding edge of (a copy of) the tile $s_{i+1}$ and get a polygon.  Since $f:P_i'\to\widetilde{Q}$ is a local embedding (the induction hypothesis), we have a local embedding from the interior of this polygon to $\widetilde{Q}$.  Then, if two adjacent edges (of tiles) in the boundary of this polygon are mapped to the same 1-cell in $\mathcal{L}$, we identify these two adjacent edges (of tiles) and get another polygon.  By eliminating such adjacent edges, we eventually get a polygon $P_{i+1}'$ and a cellular map $f: P_{i+1}'\to\widetilde{Q}$ which is a local embedding.   We denote $P_m'$ by $P^+$.  Clearly, each $V_i$ lies in the interior of $P^+$, $f:P^+\to\widetilde{Q}$ is a cellular map and a local embedding, and $P^+$ has an induced hyperbolic metric and tiling by geodesic arcs.  Note that, by this construction, it is possible that $f(s_i)$ and $f(s_j)$ are the same tile in $\widetilde{Q}$ but $s_i$ and $s_j$ are different tiles in $P^+$ (we will show later that this cannot happen).  

We first claim that $s_i$ and $s_j$ are different tiles in $P^+$ ($i\ne j$) and they do not share an edge in $P^+$.   Note that the intersection of $s_i$ with the triangular tile $t_i$ is a single point $V_i$ in $P^+$ and $t_i\ne t_j$ by our assumption.  For each $i$, we denote the edge of $t_i$ that does not contain $V_i$ by $\eta_i$.  If either $s_i=s_j$ or $s_i$ and $s_j$ share an edge in $P^+$, there is a geodesic arc $\alpha$ in $t_i\cup t_j\cup s_i\cup s_j\subset P^+$ passing through both $V_i$ and $V_j$ with two endpoints in $\eta_i$ and $\eta_j$ respectively.  Since $P$ is a convex polygon, there is another geodesic arc in $P$ connecting the two endpoints of $\alpha$, which gives a contradiction.  Thus, for each $i$, $s_i\cap\partial P'$ is the union of the two edges of $s_i$ that contain $V_i$, and $s_i\cap\partial P^+$ is the union of the edges of $s_i$ that do not contain $V_i$.

Next, we show that $P^+$ is convex. Since $P^+$ is a hyperbolic polygon with geodesic edges, it suffices to show the inner angle at each vertex of $P^+$ is less than $\pi$.  As above, $V_1,\dots,V_m$ are the only vertices where $P'$ has inner angle greater than $\pi$, and the $V_i$'s lie in the interior of $P^+$.  Let $e_i'$ and $e_i''$ be the two edges of $s_i$ incident to the $V_i$ for each $i$, and let $V_i'=\partial e_i'-V_i$ and $V_i''=\partial e_i''-V_i$ be the other endpoints of $e_i'$ and $e_i''$ respectively.  Since each tile is convex and $\partial s_i\cap int(P^+)=e_i'\cup e_i''-V_i'\cup V_i''$ in the conclusion above, for any vertex $V$ of $P^+$ that is not among the $V_i'$'s or the $V_i''$'s, the inner angle of $P^+$ at $V$ is less than $\pi$.  Thus, in order to prove $P^+$ is convex, we only need to show that the inner angles of $P^+$ at $V_i'$ and $V_i''$ (for each $i$) are less than or equal to $\pi$.  By the argument above, $V_i'$ and $V_i''$ lie in both $\partial P^+$ and $\partial P'$.  If $V_i'$ lies in the interior of an edge of $P'$, as shown in Figure~\ref{triangle}~(a), then there is a tile $\epsilon$ in $P'-int(P)$ sharing the edge $e_i'$ with $s_i$ and sharing an edge with $t_i$.  Since $V_i'$ lies in the interior of an edge of $P'$, the edge of $\epsilon$ containing $V_i'$ (other than $e_i'$) must have the other endpoint in $P$.  Since both $P$ and $\epsilon$ are convex, $\epsilon$ must be a quadrilateral tile sharing an edge with the triangular tile $t_i$, as shown in Figure~\ref{triangle}~(a), which contradicts Observation~\ref{O3}.  So, $V_i'$ must be a vertex of $P'$.  If $V_i''=V_{i+1}'$ in $\partial P'$, as shown in Figure~\ref{triangle}~(b), then there is a tile $\sigma$ in $P'-int(P)$ sharing edges with $s_i$, $s_{i+1}$, $t_i$ and $t_{i+1}$ ($\sigma\cap s_i=e_i''$ and $\sigma\cap s_{i+1}=e_{i+1}'$).  Since both $P$ and $\sigma$ are convex, $\sigma$ must be a pentagon tile sharing an edge with the triangular tiles $t_i$ and $t_{i+1}$, as shown in Figure~\ref{triangle}~(b), which contradicts Observation~\ref{O3}.  Thus, for each $i$, $V_i'$ (resp. $V_i''$) must be a vertex of $P'$ and there is only one tile from the $s_i$'s incident to $V_i'$ (resp. $V_i''$).  Since this tiling comes from geodesics in $\mathcal{L}$, $V_i'$ (resp. $V_i''$) must lie in the interior of an edge of $P^+$ for each $i$.  Therefore, $P^+$ is convex, and hence $f:P^+\to\widetilde{Q}$ is an embedding.

So, we can consider $P^+$ as a convex polygon in $\widetilde{Q}$. Next, we will analyze the points of $\mathcal{L}_0$ in the interior of an edge of $P^+$.  Since we are in the case that $M$ is not of type $\Delta(4,q,2)$ or $\Delta(3,q,2)$ and we have assumed there are triangular tiles, the triangle group for $M$ must be $\Delta(3,q,r)$ with $q, r\ge 3$ ($q$ and $r$ cannot both be $3$) and the loop $C$ in the orbifold $Q$ represents $xy^{-1}$ as shown in Figure~\ref{F31}~(a).  Suppose there is a tile $\tau$ in $P^+$ with an edge $\eta$ lying in the interior of an edge of $P^+$.  Let $A$ and $B$ be the two endpoints of $\eta$, and $\eta_A$ and $\eta_B$ be the two edges of $\tau$ (other than $\eta$) containing $A$ and $B$ respectively.  Let $A'=\partial\eta_A-A$ and $B'=\partial\eta_B-B$ be the other endpoints of $\eta_A$ and $\eta_B$ respectively.  By the argument before, either $A'\in\partial P$ or $A'$ is one of the $V_i$'s above.  We denote the edge of $P^+$ that contains $\eta$ by $\eta^+$ ($\eta\subset int(\eta^+)$).  We will show next that $A$ and $B$ are the only points in $\mathcal{L}_0\cap int(\eta^+)$.

There are several possibilities.  If $A'\ne B'$ and both $A'$ and $B'$ lie in $\partial P$, since $P$ is convex, $\tau$ must be a quadrilateral tile with one edge in $\partial P$ as shown in Figure~\ref{intersection}~(a). If $A'\ne B'$ and both $A'$ and $B'$ are among the $V_i$'s, since $P$ is convex, $\tau$ is either a hexagon as shown in Figure~\ref{intersection}~(c), or a pentagon as shown in Figure~\ref{intersection}~(d).  If $A'\in\partial P$ and $B'$ is one of the $V_i$'s, since $P$ is convex, $\tau$ must be a pentagon as shown in Figure~\ref{intersection}~(e).  However, Figure~\ref{intersection}~(d) and (e) contain pentagon tiles sharing edges with triangular tiles, which contradicts Observation~\ref{O3}.  Furthermore, if there is a hexagon tile sharing an edge with a triangular tile, since we assume the triangle group of $M$ is not Euclidean nor $\Delta(3,q,2)$, the triangle group must be $\Delta(3,q,3)$ ($q\ge 4$) and the loop $C$ represents $xy^{-1}$.  The two $q$-gons in Figure~\ref{intersection}~(c) must be among the $s_i$'s above.  Since $q\ge 4$, and by our conclusion on $s_i\cap\partial P^+$ before, $A$ and $B$ are the only points of $\mathcal{L}_0$ lying in the interior of $\eta^+$ in the case of Figure~\ref{intersection}~(c).  Hence, if there are at least 3 points of $\mathcal{L}_0$ lying in the interior of $\eta^+$, $\tau$ is either a quadrilateral tile in which case $A'\ne B'$, as shown in Figure~\ref{intersection}~(a), or a triangular tile in which case $A'=B'$, as shown in Figure~\ref{intersection}~(b).  Suppose there are at least 3 points of $\mathcal{L}_0$ lying in the interior of $\eta^+$.  Then, there must be two tiles $\tau_1$ and $\tau_2$ in $P^+$ with edges lying in $int(\eta^+)$, and $\tau_1$ and $\tau_2$ share an edge.  By the argument above, $\tau_i$ is either a triangular tile or a quadrilateral tile ($i=1,2$).  As $\tau_1$ and $\tau_2$ share an edge, by Observation~\ref{O3}, both $\tau_1$ and $\tau_2$ must be quadrilateral tiles, which is also impossible by Observation~\ref{O2}. Therefore, in any case, we have the following:

\begin{observation}\label{O4}
The polygon $P^+$ constructed above is convex, $f:P^+\to\widetilde{Q}$ is an embedding, and there are at most 2 points of $\mathcal{L}_0$ lying in the interior of any edge of $P^+$, i.e., there are at most two lines in $\mathcal{L}$ intersecting the interior of any edge of $P^+$.
\end{observation}

\begin{figure}
\begin{center}
\psfrag{l1}{$l_1$}
\psfrag{l2}{$l_2$}
\psfrag{l'}{$l'$}
\psfrag{l''}{$l''$}
\psfrag{q-gon}{$q$-gon}
\includegraphics[width=4in]{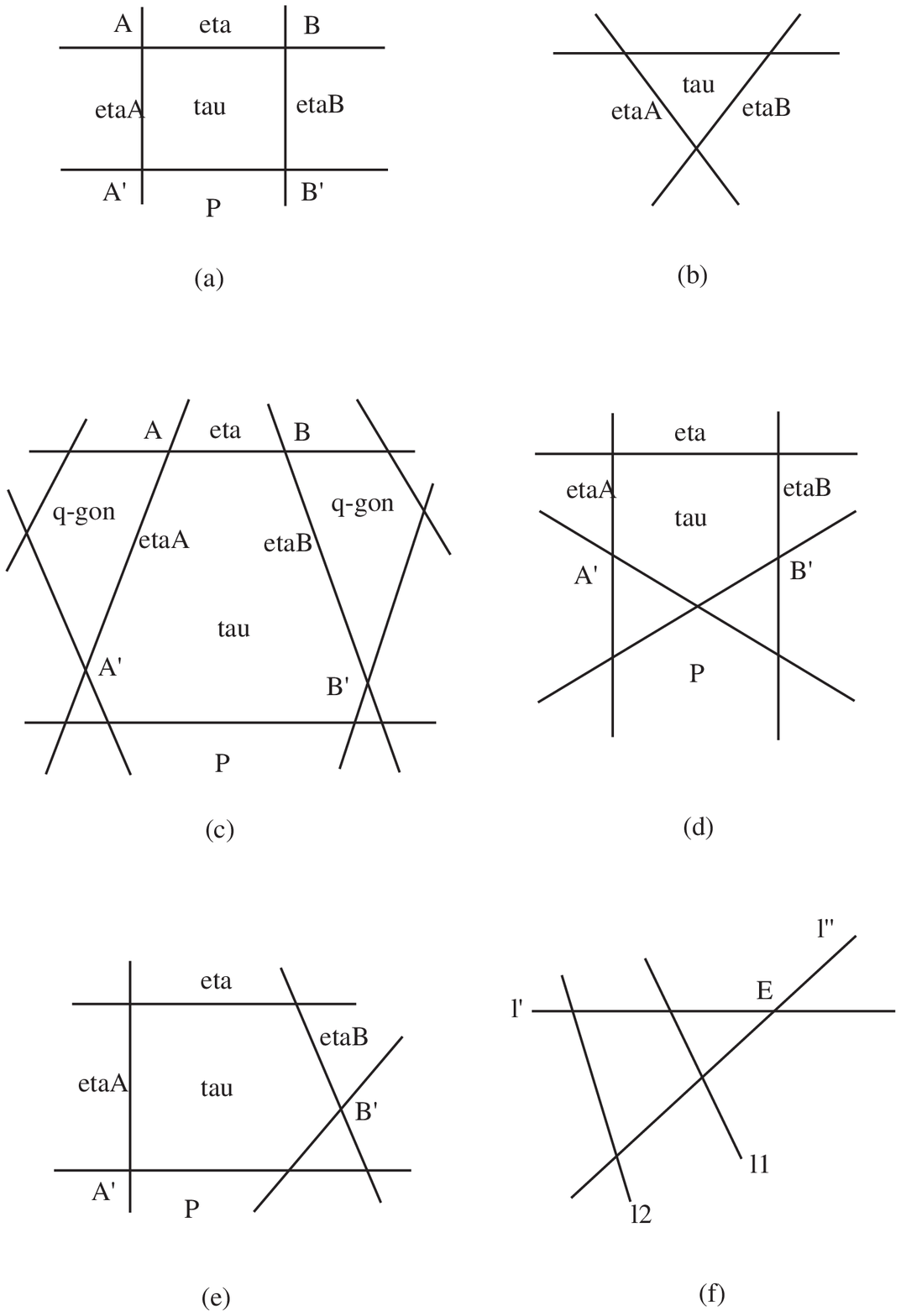}
\caption{} \label{intersection}
\end{center}
\end{figure}

Let $P_0$ be a tile in this tiling and $L_0$ be the union of lines that have nonempty intersection with $P_0$.  By Lemmas 1.7 and 1.9 of \cite{S}, each line in $L_0$ intersects at most 4 other lines in $L_0$.  Thus, if we have 5 colors to choose, we can always color a certain line in $L_0$ no matter what colors other lines (in $L_0$) have.  So, $L_0$ has the 5-color property.  We will use $P_0$ and $L_0$ as the first step of an induction to prove that $\mathcal{L}$ has the 7-color property.

\begin{subcase}\label{subcase1}  $M$ is not of type $\Delta(3,q,r)$ ($q>3$).
\end{subcase}

As before, since $M$ is not of type $\Delta(3,q,2)$ or $\Delta(4,q,2)$, lines in $\mathcal{L}$ do not pass through translates of $X$, $Y$ or $Z$ \cite{S}.  Hence, $\pi(\mathcal{L})$ is a loop as shown in Figure~\ref{F31}. 
By Lemma 1.7 in \cite{S}, there are triangular tiles in this tiling if and only if $M$ is of type $\Delta(3,q,r)$.  Thus, in this subcase, there is no triangle tile and by Observation~\ref{O2}, there are no two quadrilateral tiles sharing an edge.

The curve $C$ in this case is a figure eight, and if $M$ is of type $\Delta(p,q,r)$, the tiling consists of $p$-gons, $q$-gons and $2r$-gons.  We will start from $P_0$ and $L_0$ above.  $L_0$ has 5-color property.  Suppose we have constructed a convex polygon $P_n$. Let $L_n$ be the union of lines in $\mathcal{L}$ that have nonempty intersection with $P_n$.    Suppose $L_n$ has the 5-color property.  Let $P_{n+1}$ be the union of $P_n$ and all the tiles that have nonempty intersection with $P_n$.  Since there are no triangular tiles, by the discussion before, $P_{n+1}$ is a convex polygon embedded in $\widetilde{Q}$.  Moreover, since there are no two quadrilateral tiles sharing an edge in this tiling, by Observation~\ref{O1}, each edge of $P_{n+1}$ contains at most two points of $\mathcal{L}_0$ in its interior.

Let $L_{n+1}$ be the union of lines in $\mathcal{L}$ that have nonempty intersection with $P_{n+1}$.  Suppose we have colored $L_n$ using 5 different colors.  Let $l'\in L_{n+1}-L_n$.  Then $l'$ is a line that contains an edge $e$ of the polygon $P_{n+1}$.  Since $P_{n+1}$ is convex, $l'-P_{n+1}=l'-e$ has two components.  Next, we show that $l'-e$ does not intersect any other line in $L_{n+1}$.  Let $\gamma$ be a component of $l'-int(e)$, i.e., $\gamma$ is a ray with endpoint $E$ that is a vertex of $P_{n+1}$. Let $l''$ be the other line in $L_{n+1}$ passing through $E$.  If there is a line in $L_{n+1}$ that intersects $\gamma-E$, then it must also intersect $l''$ because $P_{n+1}$ is convex, and hence $l'$, $l''$ and this line intersect each other forming a triangle.  Since the intersection of any two lines is either empty or a single point, there must exist a triangular tile, which contradicts the hypotheses of this case.  By Observation~\ref{O1}, each edge of $P_{n+1}$ contains at most 2 points of $\mathcal{L}_0$ in its interior.  Hence, $l'$ has nonempty intersection with at most 4 other lines in $L_{n+1}$.  Since we have 5 colors to choose, we can always give $l'$ an appropriate color no matter what colors other lines (in $L_{n+1}$) have.  Therefore, $L_{n+1}$ has the 5-color property and hence $\mathcal{L}$ has the 5-color property.

\begin{subcase}\label{subcase2} $M$ is of type $\Delta(3,q,r)$, where $q, r\ge 3$.
\end{subcase}

As before, $M$ has a hyperbolic triangle group $\Delta(3,q,r)$ and $C$ is a loop representing $xy^{-1}$.  By \cite{S}, the axis $l$ of $xy^{-1}$ does not pass through any translates of $X$, $Y$ or $Z$ in this case.  Thus, $\mathcal{L}$ gives the same tiling as $\widetilde{C}$.  It is easy to see from the indices of the singular fibers that the tiling consists of triangular tiles, $q$-gon tiles and $2r$-gon tiles.  

As in subcase~\ref{subcase1}, we start with $P_0$ and $L_0$, and inductively construct a sequence of convex polygons.  Suppose we have constructed a convex polygon $P_n$.  Let $L_n$ be the union of lines in $\mathcal{L}$ that have nonempty intersection with $P_n$, and suppose $L_n$ has the 7-color property.  Let $P_n'$ be the union of $P_n$ and the tiles that intersect $P_n$.  By the discussion before, if $P_n'$ is not convex, the vertices with inner angles greater than $\pi$ come from triangular tiles.  Then, we replace $P$ and $P'$ in the argument before by $P_n$ and $P_n'$ respectively, and construct a convex polygon $P_{n+1}=P^+$ by adding some more tiles near these triangular tiles.  By Observation~\ref{O4}, the interior of any edge of $P_{n+1}=P^+$ contains at most 2 points of $\mathcal{L}_0$.

As in subcase~\ref{subcase1}, let $l'\in L_{n+1}-L_n$. Then $l'$ is a line that contains an edge $e$ of the polygon $P_{n+1}$.  Since $P_{n+1}$ is convex, $l'-P_{n+1}=l'-e$ has two components.  Let $\gamma$ be a component of $l-int(e)$.  So, $\gamma$ is a ray with one endpoint $E$ that is a vertex of $P_{n+1}$. Let $l''$ be the other line in $L_{n+1}$ passing through $E$.  If there is a line $l_1$ in $L_{n+1}$ that intersects $\gamma-E$, then it must also intersect $l''$ because $P_{n+1}$ is convex, and hence $\gamma$, $l''$ and $l_1$ intersecting each other forming a triangle.  If there is another lines $l_2$ in $L_{n+1}$ that also intersects $\gamma-E$, then $\gamma$, $l''$ and $l_2$ also intersect each other.  Since the immersed torus has 4-plane property \cite{S}, $l_1\cap l_2=\emptyset$.  So, $\gamma$, $l''$, $l_1$ and $l_2$ must form two nested triangles as shown in Figure~\ref{intersection}~(f).  However, by Lemmas 1.7 in \cite{S}, if three lines intersect each other, they must form a triangular tile, and hence such an intersection pattern is impossible.  Thus, there is at most one line in $L_{n+1}$ that intersects $\gamma-E$.  

Since there are at most 2 lines in $\mathcal{L}$ intersecting the interior of the edge $e$.  As above, there is at most one line in $L_{n+1}$ intersecting each component of $l'-e$.  Hence, there are at most 6 lines (including the two lines passing through the endpoints of $e$) in $L_{n+1}$ that intersect $l'$.  Since we have 7 colors to choose, similar to subcase~\ref{subcase1}, $L_{n+1}$ has 7-color property.

\begin{case}\label{case2}
$M$ is of type $\Delta(4,q,2)$ ($q\ge 5$).
\end{case}
For hyperbolic triangle group $\Delta(4,q,2)$ ($q\ge 5$), the axis $l$ of $xy^{-1}$ passes through translates of $X$ and $Z$.  It was shown in the section 1 of \cite{S} (see Figure 5 in \cite{S}) that the tiling induced by $\mathcal{L}$ consists only $q$-gons and the inner angle at the each vertex of a $q$-gon tile is $\pi/2$, as shown in Figure~\ref{ex}~(a).  In Figure~\ref{ex}~(a), the $X$'s, $Y$'s and $Z$'s are translates of $X$, $Y$ and $Z$ respectively.  The translates of $X$ are points in $\mathcal{L}_0$, the translates of $Y$ are centers of the $q$-gons, and the translates of $Z$ lie in $\mathcal{L}-\mathcal{L}_0$.

In fact, one can draw the lines of $\widetilde{C}$ (which is the preimage in $\widetilde{Q}$ of the figure eight loop in Figure~\ref{F31}~(a)) from the indices of the cone points.  There are parallel lines in $\widetilde{C}$, and the preimage of the cone points with indices 2 and 4 lie in the strips bounded by the parallel lines.  After collapsing every strip to a line, one gets Figure~\ref{ex}~(a).

So, there are no triangular or quadrilateral tiles in this tiling induced by $\mathcal{L}$.  Hence, by the same arguments as those in subcase~\ref{subcase1} of case~\ref{case1}, $\mathcal{L}$ (and hence $\widetilde{C}$) has the 5-color property.

\begin{figure}
\begin{center}
\psfrag{q-gon}{$q$-gon}
\includegraphics[width=4in]{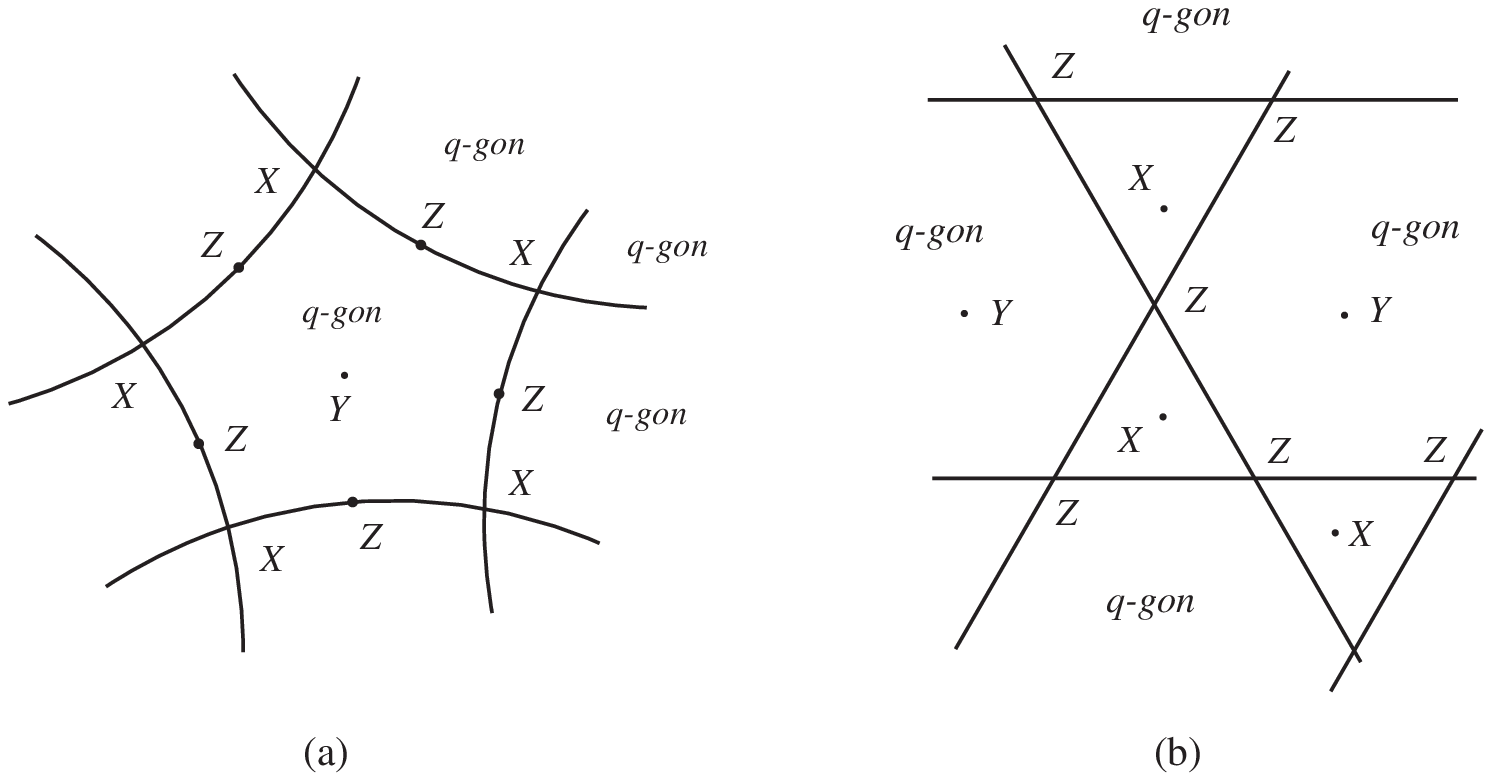}
\caption{}\label{ex}
\end{center}
\end{figure}

\begin{case}\label{case3} $M$ is of type $\Delta(3,q,2)$ ($q\ge 7$).
\end{case}

For hyperbolic triangle group $\Delta(3,q,2)$ ($q\ge 7$), the axis $l$ of $xy^{-2}$ passes through translates of $Z$, as shown in Figure 9 of \cite{S}.  The tiling induced by $\mathcal{L}$ in this case consists of $q$-gon tiles and triangular tiles, as shown in Figure~\ref{ex}~(b), and each edge of a tile is shared by a $q$-gon tile and a triangular tile.  In this tiling, as shown in Figure~\ref{ex}~(b), the translates of $Z$ are points in $\mathcal{L}_0$, the translates of $Y$ are centers of the $q$-gon tiles, and the translates of $X$ are centers of the triangular tiles.

Similar to case~\ref{case2}, one can draw the lines of $\widetilde{C}$ (which is the preimage in $\widetilde{Q}$ of the loop in Figure~\ref{F31}~(b)) from the indices of the cone points.  There are parallel lines in $\widetilde{C}$ bounding some strips in $\widetilde{Q}$.  After collapsing every strip to a line, one gets Figure~\ref{ex}~(b).  Thus, by the same arguments as those in subcase~\ref{subcase2} of case~\ref{case1}, we have that $\mathcal{L}$ (and hence $\widetilde{C}$) has the 7-color property.

Therefore, the torus $T$ has the 7-color property, where $T$ is the vertical torus corresponding to the loops in Figure~\ref{F31}.  If $M$ is of type $\Delta(3,q,2)$, $T$ has two double curves, otherwise, $T$ has only one double curve.  
\end{proof}

Let $f:S\to M$ be an immersed normal surface.  The weight of $f(S)$ is defined as $|f^{-1}(\mathcal{T}^{(1)})|$, where $\mathcal{T}^{(1)}$ is the one skeleton of the triangulation $\mathcal{T}$.  We say an immersed surface has the \emph{least weight} if the weight of $f(S)$ is minimal among all immersed surfaces in the homotopy class of $f:S\to M$.  In any homotopy class of an essential surface, there is always a least weight normal surface \cite{JR}. Weight (of a surface) is a combinatorial analogue of area.  Least area surfaces have many remarkable properties \cite{FHS}, and the results in \cite{FHS} for least area surfaces are also true for least weight surfaces \cite{JR}.  In particular, if there is an essential surface (in a certain homotopy class) having the $k$-plane and 1-line properties, then there is a least weight normal surface (in this homotopy class) that also has the $k$-plane and 1-line properties \cite{JR}.  Moreover, the intersection relation (between the planes in the universal cover) for two homotopic surfaces with the 1-line property are the same, and it follows from \cite{FHS, JR} that if an essential surface (in a certain homotopy class) has the $k$-color property, then there is a least weight normal surface (in this homotopy class) that also has the $k$-color property.  Thus, we can assume our immersed surfaces are least weight normal surfaces, and there is a least weight normal torus in the homotopy class of $T$ having the 7-color property.

\section{Immersed branched surfaces}\label{branch}

Branched surfaces have been proved fruitful in the study of incompressible surfaces and essential laminations, e.g., \cite{FO}.  For the definitions and notation related to branched surfaces, see \cite{FO, O}.  In \cite{Ch, L}, immersed branched surfaces were introduced to study immersed surfaces with small complexity.  

\begin{definition}
Let $B$ be a branched surface properly embedded in some compact 3-manifold, i.e. the local picture of $B$ in this manifold is as in Figure~\ref{F21}~(a).  Let $N(B)$ be a regular neighborhood of $B$, as shown in Figure~\ref{F21}~(b).  $N(B)$ can be considered as an $I$-bundle over $B$.  The boundary of $N(B)$ consists of vertical boundary $\partial_vN(B)$ and horizontal boundary $\partial_hN(B)$, as shown in Figure~\ref{F21} (see \cite{FO, O} for more details about branched surfaces).  Let $T$ be a surface embedded in $N(B)$.  We say that $T$ is \emph{carried} by $N(B)$ if $T$ intersects the $I$-fibers transversely.  We say  $T$ is \emph{fully carried} by $N(B)$ if $T$ transversely intersects every $I$-fiber of $N(B)$.

Let $f: B\to M$ (resp. $f: N(B)\to M$) be a map from $B$ (resp. $N(B)$) to a 3-manifold $M$.  We call $f:B\to M$ (or simply $f(B)$) an \emph{immersed branched surface} in $M$ if  $f:B\to M$ and $f:N(B)\to M$ are local embeddings.  Note that $B$ and $N(B)$ are properly embedded in some compact 3-manifold, in particular, Figure~\ref{bad} cannot be a local picture of $f(B)$.  An immersed surface $j: S\to M$ (or simply $S$) is said to be  \emph{carried} by $f:B\to M$ if, after some homotopy in $M$, $j=f\circ i$, where $i:S\to N(B)$ is an embedded surface that transversely intersects the interval fibers of $N(B)$.  We say $j: S\to M$ (or simply $S$) is \emph{fully carried} by $f:B\to M$ if $i(S)$ (as above) transversely intersects every $I$-fiber of $N(B)$.
\end{definition}

\begin{figure}
\begin{center}
\psfrag{(a)}{(a)}
\psfrag{(b)}{(b)}
\psfrag{horizontal}{$\partial_hN(B)$}
\psfrag{v}{$\partial_vN(B)$}
\includegraphics[width=4.0in]{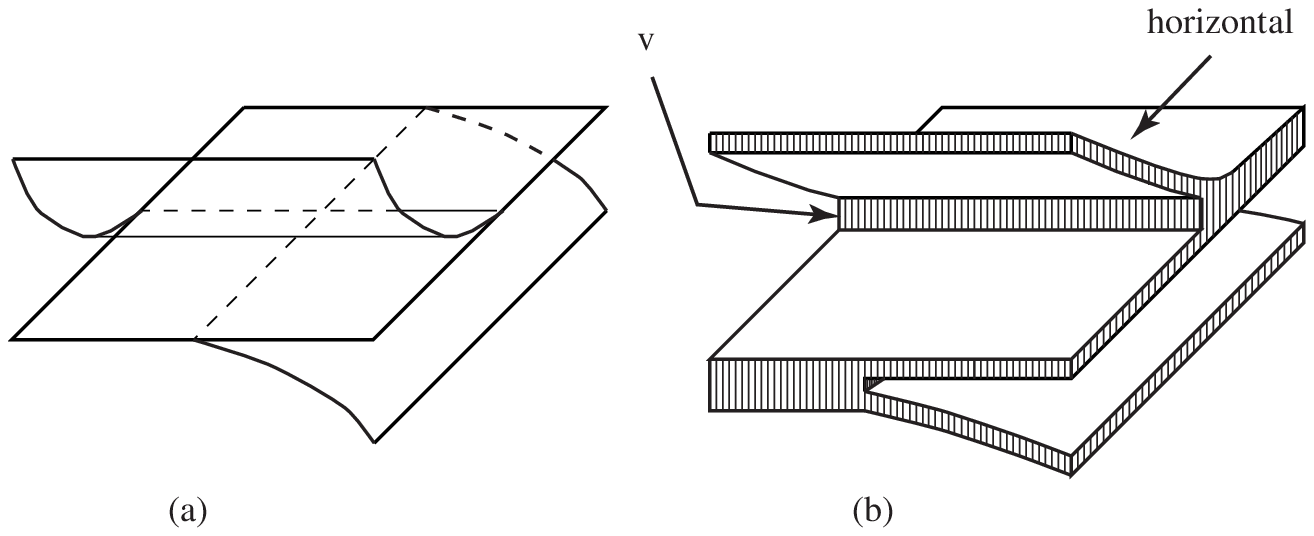}
\caption{}\label{F21}
\end{center}
\end{figure}

We assume $M$ has a certain triangulation.  Then, every $\pi_1$-injective surface can be homotoped into normal form, see \cite{JR} for some properties of normal surfaces.  In \cite{FO}, Floyd and Oertel studied embedded normal essential surfaces.  One of the results in \cite{FO} is that embedded normal surfaces are fully carried by finitely many embedded branched surfaces.  The proof of this result is straightforward.  Since there are $7$ different types of normal disks in a tetrahedron, by identifying the normal disks (in the normal surfaces) of the same type to a branch sector, one can easily construct finitely many embedded branched surfaces carrying all embedded normal surfaces.  However, this simple result is not true for immersed surfaces under our definition of immersed branched surfaces and carrying above, although every immersed essential surface can also be homotoped into normal form.  An analogue of the theorem of Floyd and Oertel for surfaces with the $4$-plane property can be found in \cite{L}.  Next, we show that this finiteness theorem can be generalized to immersed surfaces with the $n$-color property.  The proof is similar to the case of embedded surfaces in \cite{FO}.  We can identify normal disks of the same type and in the same color to a branch sector, and get only finitely many possible immersed branched surfaces.

\begin{lemma}\label{L31}
The immersed $\pi_1$-injective surfaces (in $M$) with the $n$-color property are fully carried by finitely many immersed branched surfaces.
\end{lemma}
\begin{proof}
Let $\mathcal{T}$ be a triangulation of $M$, $T$ be a tetrahedron in $\mathcal{T}$, and $T'$ be a tetrahedron in $\widetilde{\mathcal{T}}$ such that $\pi (T')=T$, where $\widetilde{\mathcal{T}}$ is the induced triangulation of the universal cover $\widetilde{M}$ and $\pi:\widetilde{M}\to M$ is the covering map.  Let $f:S\to M$ be a normal $\pi_1$-injective surface with the $n$-color property.  After homotopy, we can assume $f$ is a normal surface and has least weight (or combinatorial area as in \cite{JR}).  Each tetrahedron has 7 normal disk types.  Let $d$ be a normal disk type in $T'$ and $\gamma_d$ be the collection of normal disks in $\pi^{-1}(f(S))\cap T'$ of type $d$.  Suppose we have colored all the planes in $\pi^{-1}(f(S))$ using $n$ different colors.  Then, we give every normal disk in $\pi (\gamma_d)$ a color induced from $\gamma_d$.  We still call the correspondent normal disk type in $T$ the disk type $d$, and we use $d\times I$ to denote the product of an interval $I$ and a normal disk of type $d$.  As in \cite{FO, L}, we put $k$ ($k\le n$ is the number of different colors in $\gamma_d$) such products $d\times I$'s in $T$ such that any two disks in the same color are in the same $d\times I$ and are transverse to the $I$-fibers.  We can do such construction for every normal disk type in every tetrahedron in $T$.

Since our coloring in $\widetilde{M}$ is not equivariant, we must be careful when we glue these $d\times I$'s together along the 2-skeleton (to form a fibered neighborhood of an immersed branched surface). Let $N_k$ be a small neighborhood of the $k$-skeleton and $N$ be a small neighborhood of a 2-simplex in $M-N_1$.  There are finitely many $d\times I$'s intersecting $N$.  We say that two disks in $f(S)\cap N$ are of the same type if they connect the same two products (i.e. the same two $d\times I$'s) in the two tetrahedra sharing this 2-simplex.  Thus, we can put finitely many product regions $E_2\times I$'s in $N$ such that any two disks in $f(S)\cap N$ of the same type are in the same product region.  We can do such construction for every 2-simplex.  Similarly, we can put some product regions in $N_1$ connecting those product regions in $M-N_2$ and $N_2-N_1$.  Note that the number of different product regions for a 1-simplex depends on $n$ and the number of tetrahedra incident to this 1-simplex.  As in \cite{FO}, by identify every $I$-fiber of every product region above to a point, we can construct a singular branched  surface and the union of those product regions can be consider as a fibered neighborhood of this singular branched surface.  By our construction above, two normal disks in the immersed surface are identified to the same branch sector only if they have the same induced color, and hence the two corresponding planes in the universal cover have the same color and do not intersect each other.  Thus, there is no local picture of Figure~\ref{bad} in this singular branched surface.  By \cite{C}, this singular branched surface can be embedded in some 3-manifold, i.e., one can construct a 3-manifold $M'$ in which $B$ is embedded, and there is a local embedding $h:B\to M$ such that $h(B)$ is the singular branched surface we constructed.  Moreover, since the number of those product regions are bounded, there are only finitely many ways to construct such immersed branched surfaces, and the lemma follows.
\end{proof}

\begin{figure}
\begin{center}
\resizebox{3in}{!}{\includegraphics{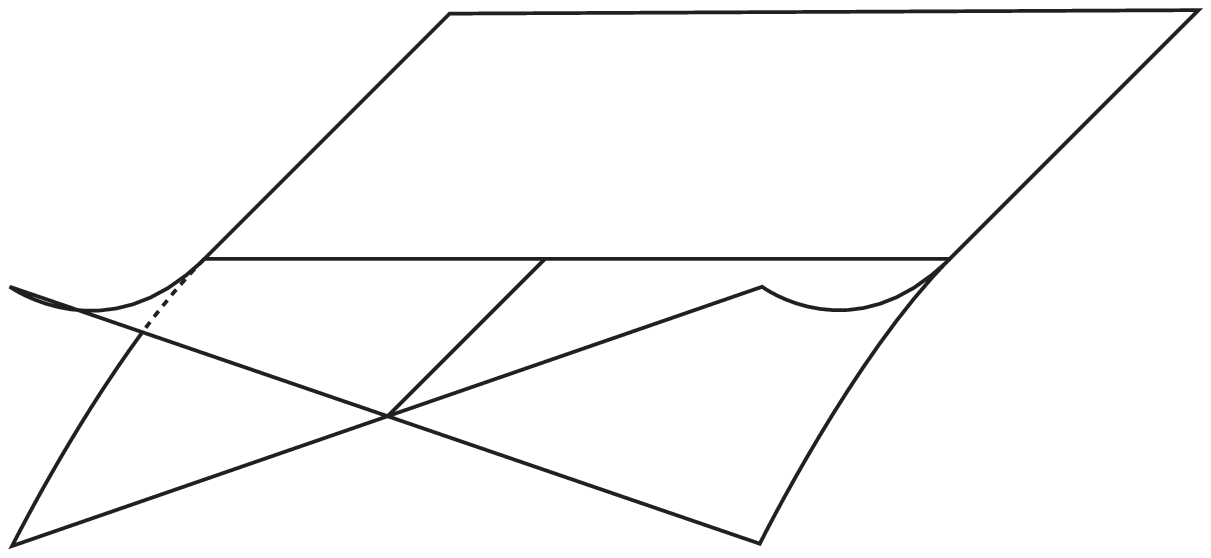}}
\caption{} \label{bad}
\end{center}
\end{figure}

In the proof of Lemma~\ref{L31}, the number of possible immersed branched surfaces depends on $n$ and the triangulation.  

Let $g:S\to M$ be a least weight $\pi_1$-injective surface and $f: N(B)\to M$ be (a fibered neighborhood of) the immersed branched surfaces fully carrying $S$ that we constructed in the proof of Lemma~\ref{L31}.  By our definition of carrying, we can view $S$ as an embedded surface in $N(B)$ and $f|_S=g$.  Let $D\times I$, where $D$ is a disk, be a product region in $N(B)$ with each $I$-fiber a subarc of an $I$-fiber of $N(B)$, $D\times\partial I\subset S$ and $S\cap (D\times int(I))=\emptyset$.  Note that by our construction, $f(D\times int(I))\cap f(S)$ may not be empty.  Nevertheless, the following two properties follow trivially from our construction in the proof of Lemma~\ref{L31}.

\begin{property}\label{P1}
Let $f$, $N(B)$, and $D\times I$ be as above, and suppose $f|_{D\times I}$ is an embedding.  Let $\tilde{f}: D\times I\to\widetilde{M}$ be a lift of $f|_{D\times I}$, and we denote the two planes containing the two components of $\tilde{f}(D\times\partial I)$ by $S_1$ and $S_2$.  Then, by our construction, we may assume $S_1$ and $S_2$ have the same color, and hence either $S_1=S_2$ or $S_1\cap S_2=\emptyset$.
\end{property}

\begin{property}\label{P2}
Let $S_i$ ($i=1,2$) and $\tilde{f}$ be as in Property~\ref{P1}.  So, we may assume $S_1$ and $S_2$ have the same color.  If there is a plane $S_3$ that intersects $\tilde{f}(D\times int(I))$, then $S_3$ cannot be in the same color as $S_1$ (since we assumed $S\cap D\times int(I)=\emptyset$ in $N(B)$ before) and in particular, $S_3\ne S_1$. 
\end{property}

Floyd and Oertel also showed \cite{FO} that the branched surface fully carrying a least weight embedded incompressible surface can be split to be an incompressible branched surface.  In particular, it contains no monogon.  This result can be generalized to immersed surfaces as well.

\begin{definition}
An immersed branched surface $f:B\to M$ (and $f:N(B)\to M$) is said to be \emph{incompressible} if it satisfies the following conditions.
\begin{enumerate}
\item  $f|_{\partial_hN(B)}$ is $\pi_1$-injective.

\item $f(B)$ has no immersed monogon, i.e. for any closed curve $c=\alpha\cup\beta$, where $\alpha$ is a vertical arc in $\partial_vN(B)$ and $\beta\subset\partial_hN(B)$, $f|_c$ is homotopically nontrivial.
\end{enumerate}
\end{definition}

Note that, if $f|_{\partial_hN(B)}$ is $\pi_1$-injective, then $B$ dose not contain any disk of contact (a disk of contact \cite{FO} is a disk $D$ embedded in $N(B)$, transverse to the $I$-fibers of $N(B)$, and $\partial D\subset\partial_vN(B)$).  A disk of contact can be eliminated by splitting $B$ \cite{FO}.  One can easily modify an argument in \cite{FO} to prove the following lemma. 

\begin{lemma}\label{L22} Suppose $M$ is a closed orientable 3-manifold with trivial $\pi_2(M)$.  Let $S$ be a closed orientable $\pi_1$-injective least weight normal surface.  Suppose $S$ is fully carried by an immersed normal branched surface $f:B\to M$.  Then, after eliminating disks of contact in $B$ and taking a sub branched surface if necessary, $f:B\to M$ is incompressible.
\end{lemma}
\begin{proof}
Since $S$ is fully carried by $f:B\to M$, we can assume $S$ lies in $N(B)$.  After some isotopy and taking multiple copies of $S$ if necessary, we can also assume $\partial_hN(B)\subset S$.  Since $S$ is $\pi_1$-injective, after getting rid of disks of contact (by splitting $B$) and taking a sub branched surface if necessary, $f|_{\partial_hN(B)}$ is $\pi_1$-injective.  So, it suffices to show that there is no immersed monogon.   Suppose that there is a monogon, i.e., a closed curve $c=\alpha\cup\beta$ such that $f|_c$ is homotopically trivial, where $\alpha$ is a vertical arc in $\partial_vN(B)$ and $\beta\subset\partial_hN(B)\subset S$.  Let $A$ be the component of $\partial_vN(B)$ that contains $\alpha$ ($A$ is an annulus), and let $S_A$ be the horizontal boundary component that contains $\beta$ and $\partial A$.  Let $R$ be a small rectangular neighborhood of $\alpha$ in $A$, where $\partial R$ consists of two arcs in $\partial A$ and two vertical arcs $\alpha_1$ and $\alpha_2$ of the annulus $A$.  Let $\beta_i$  be an arc in $S_A$ that is parallel and close to $\beta$ such that $\partial\beta_i=\partial\alpha_i$ ($i=1,2$).  So, $\gamma =\beta_1\cup (\partial A-R)\cup\beta_2$ is a closed curve in $S$.  Since $f|_c$ is homotopically trivial, $f|_\gamma$ is homotopically trivial.  Since $f|_{S_A}$ is $\pi_1$-injective, $\gamma$ must be a trivial curve in $S_A$.  Since there is no disk of contact and $f|_{\partial_hN(B)}$ is $\pi_1$-injective, $S_A$ must be an annulus.  Moreover, $f(S_A)$ can be homotoped into $f(A)$ fixing $f(\partial A)$, since $\pi_2(M)$ is trivial.  After this homotopy, the total weight of $f(S)$ is reduced, which contradicts the assumption that $S$ has least weight.
\end{proof}

Note that all disks of contact can be found algorithmically by solving a system of branch equations (see \cite{FO,O} for the definition of branch equation and \cite{AL, JO} for more details).  In fact, one can formulate a system of nonhomogeneous linear equations for surfaces (carried by the branched surface) with boundary in the branch locus and having Euler characteristic 1.  By solving such linear systems and finding solutions corresponding to disks of contact with least weight \cite{AL,JO}, one can successively find all circles in the branch locus that bound disks of contact.  After successively splitting along these disks of contact, we get a branched surface without disk of contact and still carrying the surface $S$.  After taking a sub branched surface if necessary, we get a branched surface fully carrying $S$. 

Therefore, we can algorithmically construct finitely many immersed incompressible branched surfaces such that each least weight $\pi_1$-injective normal surface with the 7-color property is fully carried by one of them.  

\section{Essential tori}\label{tori}
Suppose $B$ is an embedded branched surface in a 3-manifold $N$, and $f:B\to M$ ($f:N(B)\to M$) be an immersed branched surface.  Suppose there is a least weight normal essential torus $T$ fully carried by $f:B\to M$. In other words, $T\subset N(B)\subset N$ is a torus embedded in $N(B)$ and transversely intersecting every $I$-fiber of $N(B)$, and $f|_T$ is an immersed essential normal torus that has least weight in its homotopy class.

As $f(T)$ is a normal surface, $f^{-1}(\mathcal{T}^{(2)})$ gives a cell decomposition of $N(B)$ and the weight of $f(T)$ is equal to the number of intersection points of $T\cap f^{-1}(\mathcal{T}^{(1)})$.  We denote the weight of $f(T)$ by $weight(T)$.  We call an immersed essential torus $T$ in $M$ an \emph{absolutely least weight essential torus} if the weight of any essential torus in $M$ is greater than or equal to $weight(T)$.

\begin{lemma}\label{bound}
Let $M$ be an orientable small Seifert fiber space of hyperbolic type.  Suppose  $f:B\to M$ is an immersed branched surface (constructed in section~\ref{branch}) fully carrying a least weight normal essential torus with the $n$-color property.  Then, there is a number $W$ that depends only on $f:B\to M$ and can be algorithmically calculated, such that $M$ contains some essential torus with weight less than $W$, i.e., $W$ is an upper bound on the weight of an absolutely least weight torus.
\end{lemma}
\begin{proof}
Let $T$ be the least weight normal essential torus with the $n$-color property and suppose $f:B\to M$ is constructed in section~\ref{branch} using $T$.  By Lemma~\ref{L22}, we can assume that $f:B\to M$ is incompressible.  Suppose $N(B)$ is embedded in a closed 3-manifold $N$, and $T$ is embedded in $N(B)$ transverse to every $I$-fiber of $N(B)$.  Since $f|_{\partial_hN(B)}$ is $\pi_1$-injective, similar to embedded branched surfaces, we can assume every 2-sphere component of $\partial N(B)$ bounds a 3-ball (in $N$) which is in the form $D^2\times I$, where $D^2\times\partial I$ consists of two disk components of $\partial_hN(B)$ and $\partial D^2\times I$ is a component of $\partial_vN(B)$.  Since $M$ is irreducible, we can extend the map $f$ through these $D^2\times I$ regions.

\begin{S4case}\label{C1}
Every boundary component of $N(B)$ is a 2-sphere.
\end{S4case}
As above, we can assume that $N$ is the union of $N(B)$ and the $D^2\times I$ regions.  Since $T$ is two-sided in $M$, $T$ is two-sided in $N$.  The $I$-fibers of $N(B)$ and the $I$-fibers of these $D^2\times I$ regions are connected together forming a one-dimensional foliation of $N$.  Moreover, since $N(B)$ fully carries $T$, $\overline{N-T}$ (i.e. the closure of $N-T$ under the path matric) is an $I$-bundle over a closed (possibly not connected) surface.  Since $T$ is two-sided, if $T$ is nonseparating in $N$, then $\overline{N-T}$ is of the form $T\times I$, and if $T$ is separating, $\overline{N-T}$ consists of two twisted $I$-bundles over Klein bottles.

We first consider the case that $T$ is nonseparating in $N$, i.e., $\overline{N-T}=T\times I$.  Since $N-N(B)$ is a union of 3-balls and $M$ is irreducible, the map $f:N(B)\to M$ can be extended to a map $f:N\to M$.  Next, we show that $f:N\to M$ is $\pi_1$-injective.  Otherwise, there is a closed essential curve $\gamma$ in $N$ such that $f(\gamma)$ is homotopically trivial.  After some homotopy, we can assume $\gamma$ is a union of  $2k$ arcs $\alpha_1, \beta_1,\dots,\alpha_k,\beta_k$ such that each $\alpha_i$ lies in $T$ and each $\beta_i$ is an $I$-fiber of $\overline{N-T}=T\times I$.  Moreover, by fixing a direction for $\gamma$ and a normal direction for $T$, we can assume (after homotopy) that the induced direction (from the fixed direction of $\gamma$) for each $\beta_i$ agrees with the normal direction of $T$.  We call these $\beta_i$'s $\beta$-arcs.  As $f|_T$ is $\pi_1$-injective, there is at least one $\beta$-arc in $\gamma$.  Since $T$ is fully carried by $N(B)$, i.e. $T$ is embedded in $N(B)$ and transversely intersects every $I$-fiber of $N(B)$, each $I$-fiber of $\overline{N-T}=T\times I$ is a union of subarcs of $I$-fibers in $N(B)$ and at most one $I$-fiber of some $D^2\times I$ region in $N-int(N(B))$.  Since each $I$-fiber in a $D^2\times I$ region can be isotoped to a subarc of an $I$-fiber of $N(B)$, after some homotopy, we can assume each $\beta_i$ is a subarc of an $I$-fiber of $N(B)$.

Let $\gamma'$ be a lift of $f(\gamma)$ in the universal cover of $M$.  Since $f:N(B)\to M$ is a local embedding and $f(\gamma)$ is homotopically trivial, $\gamma'$ is a union of $2k$ arcs $\alpha_1',\beta_1',\dots,\alpha_k',\beta_k'$, where the $\alpha_i'$'s and $\beta_i'$'s are corresponding lifts of the $\alpha_i$'s and $\beta_i$'s.  Thus, each $\alpha_i'$ lies in a plane in the preimage of $f(T)$ and each $\beta_i'$ is an $I$-fiber of the preimage of $f(N(B))$.  Moreover, each plane in the preimage of $f(T)$ has a normal direction induced from that of $T$, and the induced direction for each $\beta_i'$ is compactible with the induced normal directions of these planes.  

For each $\beta_i'$, we use $T_i^+$ and $T_i^-$ to denote the two planes in $\pi^{-1}(f(T))$ containing the two endpoints of $\beta_i'$ respectively with $T_i^+=T_{i+1}^-$ for each $i$ ($T_{k+1}^{\pm}= T_1^{\pm}$). Note that $T_i^+$ and $T_i^-$ cannot be the same plane, otherwise the directions will be incompatible.  Moreover, by Property~\ref{P2} in the last section, $T_i^{\pm}\cap int(\beta_i')=\emptyset$.  By Property~\ref{P1}, $T_i^+\cap T_i^-=\emptyset$.  Since each plane in $\pi^{-1}(f(T))$ is embedded and separating in $\widetilde{M}$, these $T_i^{\pm}$'s must be mutually disjoint planes that cut $\widetilde{M}$ into $k+2$ pieces.  Since the induced direction of $\gamma'$ is compatible with the induced normal direction of these $k+1$ planes, and since $T_i^{\pm}\cap int(\beta_i')=\emptyset$, these $\beta_i'$'s must lie in different pieces, and this is impossible as $\gamma'$ is a closed curve in $\widetilde{M}$

Therefore, in the case that $T$ is nonseparating, $f:N\to M$ is $\pi_1$-injective and in particular $f|_H$ is $\pi_1$-injective for any torus $H$ carried by $N(B)$.  By solving the system of branched equations (see \cite{AL,JO} for more details), one can always find a solution corresponding to a torus fully carried by $N(B)$, and its weight is an upper bound on the weight an absolutely least weight torus.  Note that, by solving the system of branch equations, one can always find a certain torus fully carried by a branched surface, but one may not be able to find the essential torus $T$ that we use to construct the immersed branched surface in section~\ref{branch}.

Now, we suppose $T$ is separating, i.e., $\overline{N-T}$ consists of two twisted $I$-bundles over Klein bottles.  In this case, $N$ has a double cover $\hat{N}$ where $T$ lifts to a nonseparating torus.  Let $\hat{T}$ and $N(\hat{B})$ be the preimage of $T$ and $N(B)$ in this double cover, and $g:\hat{N}\to M$ be the composition of the map $f$ and this double covering.  Then, similar to the argument above, $g$ induces an injection on $\pi_1$.  As above, one can solve the system of branch equations and find a torus $H$ carried by $N(\hat{B})$.  The map $g$ restricted to $H$ is $\pi_1$-injective, and the weight of $H$ is an upper bound on the weight of an absolutely least weight torus.

\begin{S4case}
Some boundary component of $N(B)$ is not a 2-sphere.
\end{S4case}
Since $f|_{\partial_hN(B)}$ is $\pi_1$-injective, there must be a component $A$ of $\partial_hN(B)$ that is not a disk.  Let $X$ be an arbitrary torus fully carried by $N(B)$.  Note that we can find such a torus $X$ by solving the system of branch equations.  After some isotopy on $X$, we can assume $A\subset X$.  Since $N(B)$ has no disk of contact, $\partial A$ must consist of essential curves in $X$.  Hence, $A$ must be an annulus and every nondisk component of $\partial_hN(B)$ is an annulus.  Let $X'$ be a disjoint union of some parallel copies of $X$.  We can choose $X'$ so that each nondisk component of $\partial_hN(B)$ lies in some component of $X'$ and each component of $X'$ contains some nondisk component of $\partial_hN(B)$.  Let $N_B$ be the union of $N(B)$ and all the $D^2\times I$ regions of $N-int(N(B))$.  By our assumptions on $X'$, $X'-\partial N_B$ consists of annuli.  We can cut $N_B$ open along $X'$, and by the discussion above, we get a union of $I$-bundles over some compact surfaces with Euler characteristic zero.  Let $Y$ be another torus fully carried by $N(B)$.  After some isotopy in $N_B$, we can assume that the number of intersection curves of $Y$ and $X'$ is minimal.  Since both $X$ and $Y$ are $\pi_1$-injective in $N(B)$, $X\cap Y$ must consist of essential curves in both $Y$ and $X'$.  Moreover, since $N_B$ is cut by $X'$ into a union of $I$-bundles (over some compact surfaces with Euler characteristic zero), the intersection of $Y$ with each $I$-bundle must consist of annuli that can be homotoped either into $X'$ or to a vertical annulus in this $I$-bundle.  As $X'-\partial N_B$ consists of annuli and $Y\cap X'\subset X'-\partial N_B$, we can choose some disjoint vertical annuli in $\overline{N_B-X'}$ so that $Y$ can be homotoped to a possibly singular torus lying in the union of these vertical annuli and $X'$.  The union of $X'$ and these disjoint vertical annuli form a 2-complex $C_Y$.  Since these vertical annuli are disjoint and $Y\cap X'\subset X'-\partial N_B$, the 1-skeleton (i.e. the boundary of these vertical annuli) of this 2-complex $C_Y$ consists of disjoint circles in $X'-\partial N_B$.  In fact, by choosing enough such disjoint vertical annuli in $\overline{N_B-X'}$, we can construct a 2-complex $C$ in $N(B)$, which is a union of $X'$ and some disjoint vertical annuli in $\overline{N_B-X'}$, such that every torus fully carried by $N(B)$ can be homotoped to a possibly singular torus in this 2-complex $C$ and the 1-skeleton (i.e. the boundary of these vertical annuli) of $C$ consists of disjoint circles in $X'-\partial N_B$.  Clearly, $C$ can be found algorithmically.

Let $A$ be an annular component of $\partial_hN(B)$.  Since $X'-\partial N_B$ is a union of annuli, the boundaries of those vertical annuli in the construction of $C$ above consist of circles (in $X'$) parallel to $\partial A$.  By mapping circles parallel to $\partial A$ to points, we can project $X'$ to a union of circles and project those vertical annuli in $C$ to arcs connecting these circles.  So, we get a graph $G$ such that $C=G\times S^1$.  By the discussion above, for each torus $Y$ fully carried by $N(B)$, there is a closed curve $\gamma_Y$ in $G$ such that $\gamma_Y\times S^1\subset G\times S^1\subset N(B)$ is homotopic to $Y$ in $N_B$.  As $G\times S^1\subset N(B)$, we have an induced map $f:G\times S^1\to M$.  By our hypotheses, there is a torus $T$ fully carried by $N(B)$ and $f|_T$ is $\pi_1$-injective.  Hence, we know that there is a closed curve $\eta$ in $G$ such that $f:\eta\times S^1\to M$ is $\pi_1$-injective.  

We can express each element in $\pi_1(M)$ in the form of $al^k$, where $a\in\Delta(p,q,r)$ and $l$ is a generator of the cyclic center of $\pi_1(M)$.  Let $(x,y)$ ($x\in G, y\in S^1$) be a base point in $G\times S^1$. Suppose $f(\{x\}\times S^1)$ represents an element $al^k$ and $f(\eta\times\{y\})$ represents $bl^m$, where $a,b\in\Delta(p,q,r)$.  Thus,  $al^k$ and $bl^m$ commute in $\pi_1(M)$, which implies that $a$ and $b$ must commute in the hyperbolic triangle group $\Delta(p,q,r)$, and hence $a$ and $b$ must generate a cyclic subgroup in $\Delta(p,q,r)$. 

\begin{S4subcase}\label{sb1}
$a\ne 1$.
\end{S4subcase}

In this case, since $<a,b>$ is a cyclic subgroup in $\Delta(p,q,r)$ and since $f:\eta\times S^1\to M$ is $\pi_1$-injective, $a$ must have infinite order in $\Delta(p,q,r)$. 

Since $G$ is a graph, $\pi_1(G)$ is a free group generated by $g_1,\dots, g_n$.  Let $\gamma_1,\dots,\gamma_n$ be closed curves in $G$ representing $g_1,\dots, g_n$ respectively.  Suppose $f(\gamma_i\times\{y\})$ represents the element $c_il^{s_i}$ in $\pi_1(M)$ ($i=1,\dots,n$), where $c_i\in\Delta(p,q,r)$.  By our construction, each $c_il^{s_i}$ commutes with $al^k$, and hence each $c_i$ commute with $a$ in the hyperbolic triangle group $\Delta(p,q,r)$.  Thus, $<c_i,a>$ is a  cyclic subgroup of $\Delta(p,q,r)$ for each $i$.  Since $a$ has infinite order in $\Delta(p,q,r)$, $<c_1,\dots,c_n,a>$ must also be an infinite cyclic subgroup in $\Delta(p,q,r)$.  If $f|_{\gamma_i\times S^1}$ is not $\pi_1$-injective for each $i$, then $<c_il^{s_i},al^k>$ is infinite cyclic in $\pi_1(M)$ for each $i$, and $<c_1l^{s_1},\dots,c_nl^{s_n}, al^k>$ must be a cyclic subgroup of $\pi_1(M)$, i.e., $f_*(\pi_1(G\times S^1))$ is cyclic in $\pi_1(M)$, which contradicts the hypotheses that there is a curve $\eta\subset G$ such that $f|_{\eta\times S^1}$ is $\pi_1$-injective.  Therefore, for some $1\le j\le n$, $f|_{\gamma_j\times S^1}$ is $\pi_1$-injective.  So, the maximal weight of the tori $f(\gamma_1\times S^1),\dots, f(\gamma_n\times S^1)$ is an upper bound on the weight of an absolutely least weight $\pi_1$-injective torus.  Since the $\gamma_i$'s can be found easily from the graph $G$, this bound can be calculated algorithmically.

\begin{S4subcase}
$a=1$.
\end{S4subcase}
In this case, $f(\{x\}\times S^1)$ represents $l^k$ and $k\ne 0$.  As $f|_{\eta\times S^1}$ is $\pi_1$-injective and $f(\eta\times\{y\})$ represents $bl^m$, $b$ must have infinite order in $\Delta(p,q,r)$.  Let $\gamma_1,\dots,\gamma_n$ and $c_1,\dots,c_n$ be as in subcase~\ref{sb1} above.  If $c_j$ has infinite order for some $j$, then $f|_{\gamma_j\times S^1}$ is $\pi_1$-injective and we get a bound as in subcase 1 above.  So, we assume each $c_i$ has finite order in $\Delta(p,q,r)$.

If $c_i$ and $c_j$ generate a cyclic group in $\Delta(p,q,r)$ for each pair $i,j$, then $<c_1,\dots,c_n>$ is a finite cyclic subgroup of $\Delta(p,q,r)$, which contradicts the assumption that $f|_{\eta\times S^1}$ is $\pi_1$-injective.  Thus, there must be a pair $c_i$, $c_j$ ($1\le i<j\le n$) such that $c_i$ and $c_j$ both have finite order but $<c_i,c_j>$ is an infinite group.  Then, by Lemma~\ref{Scott2}, there is an infinite-order element $w$ in $<c_i,c_j>$ that is of the form of $u$, or $uv^{\pm 1}$, or $uv^{\pm 2}$, where $u$ and $v$ are among $c_i$, $c_j$ and $c_ic_j$.   Let $\gamma_w$ be a loop in $G$ representing $w$.  Hence, $f|_{\gamma_w\times S^1}$ is $\pi_1$-injective.  One can enumerate all possible loops for $\gamma_w$ using the $c_i$'s above and get an upper bound on the weight of $f(\gamma_w\times S^1)$, which is also an upper bound on the weight of an absolutely least weight essential torus.
\end{proof}

\section{An algorithm}\label{S5}
In this section, we summerize our algorithm to recognize Seifert fiber spaces.  Since there are algorithms to decide whether a 3-manifold is reducible or Haken and whether a Haken manifold is a Seifert fiber space \cite{JT}, we can assume our manifold $M$ is irreducible and non-Haken.  There are algorithms \cite{R, Th} to decide whether a 3-manifold is a 3-sphere.  Rubinstein has also given an algorithm to find strongly irreducible genus $g$ Heegaard splittings, and using that he gave an algorithm to recognize lens spaces \cite{R}.  Using Rubinstein's techniques and genus 2 strongly irreducible Heegaard splitting, one can find an algorithm to decide whether a 3-manifold is a Seifert fiber space with finite fundamental group (see also \cite{RR}).  

In fact, if the conjecture is true that finite group actions on lens spaces are standard, then one has a trivial algorithm to recognize Seifert fiber spaces with finite fundamental group as follows.  If $M$ is an orientable Seifert fiber space with finite fundamental group, then there is a $k$-fold ($k\le 60$) cover of $M$ that is a lens space \cite{Or}.  So, one only needs to enumerate all possible $60$-fold covers of a 3-manifold and check each one of them to see if it is a lens space using the algorithm in \cite{R}.  

Therefore, it remains to find an algorithm to recognize non-Haken small Seifert fiber spaces with infinite fundamental groups.  

In the following algorithm to recognize small Seifert fiber spaces with infinite fundamental group, we use immersed branched surfaces discussed previously.  Since one can describe a branched surface by a system of branched equations as in \cite{AL, FO, O}, many geometric operations in the algorithm can be described by linear equations, and this algorithm can be much more simplified during implementation.

\begin{step}\label{step2}
Check whether $M$ is a small Seifert fiber space of Euclidean type.
\end{step}

If $M$ is a small Seifert fiber space of Euclidean type, then the triangle group of the corresponding orbifold is either $\Delta(2,4,4)$, or $\Delta(2,3,6)$ or $\Delta(3,3,3)$.  If $M$ is of type $\Delta(2,4,4)$ (resp. $\Delta(2,3,6)$, then $M$ is double (resp. triple) covered by a Haken Seifert fiber space whose orbifold is a 2-sphere with 4 cone points.  If $M$ is of type $\Delta(3,3,3)$, then $M$ is a triple covered by a Haken Seifert fiber space whose orbifold is a torus without cone point.  We can enumerate all possible double and triple covers of $M$.  Then, using the algorithm in \cite{JT}, we check these covers to see whether one of them is a Haken Seifert fiber space, and we will know immediately whether $M$ is a Seifert fiber space.  In fact, one can find a least weight embedded incompressible vertical torus in each of such covers using \cite{JT}.  The projection of this embedded torus to $M$ is a $\pi_1$-injective torus with the 4-plane and 1-line properties \cite{S, FHS}.  Using the 4-plane and 1-line properties, as in \cite{S, HS} (and Step~\ref{step5} below), one can algorithmically perform some simple homotopies to remove all triple points of this immersed essential torus.  By checking the complement of this tori without triple points, one can easily see the Seifert fiber structure on $M$.

\begin{step}\label{step3}  
Construct finitely many immersed branched surfaces so that any least weight essential surface with the 7-color property is fully carried by one of them.
\end{step}

For each normal disk type $\delta$, we put $k$ ($k\le 7$) normal disks of type $\delta$.  Then, as in the proof of Lemma~\ref{L31}, we can split each normal disk near the 2-skeleton and glue them together to form an immersed branched surface.  Since the number of normal disks are bounded, we can enumerate all possible immersed branched surfaces.  

\begin{step}\label{step4}  Calculate a number $W$ such that, if $M$ is a small Seifert fiber space of hyperbolic type, then $M$ contains an essential torus with weight less than $W$. 
\end{step}

By Theorem~\ref{T31}, $M$ contains an essential torus with the $7$-color property if $M$ is a small Seifert fiber space of hyperbolic type.  So, one of the immersed branched surfaces in Step~\ref{step3} fully carries a least weight normal essential torus.  By solving systems of branch equations, we can find all disks of contact in these branched surfaces. After splitting along the disks of contact and taking a sub branched surface if necessary, by Lemma~\ref{L22}, we can assume one of the immersed branched surfaces is an incompressible immersed branched surface fully carrying a least weight normal essential torus.  Therefore, for each immersed branched surface in Step~\ref{step3}, we can algorithmically calculate a number (as in the proof of Lemma~\ref{bound}) that is an upper bound on the weight of some essential torus, if this immersed branched surface fully carries an essential torus.  By taking the maximum, we obtain a number $W$, which is an upper bound on the weight of an absolutely least weight essential torus if the 3-manifold is a small Seifert fiber space of hyperbolic type.

\begin{step}\label{step5}
Determine whether $M$ is a small Seifert fiber space of hyperbolic type.
\end{step}

Let $W$ be the number calculated in Step~\ref{step4}.  Then, by a theorem in \cite{LN} (i.e., Theorem~\ref{T:LN} as follows), in a Seifert fiber space of hyperbolic type, there is always a torus $T$ with the 4-plane and 1-line properties and $weight(T)\le W$.

\begin{theorem}[\cite{LN}]\label{T:LN} In a non-Haken small Seifert fiber space with infinite $\pi_1$, an absolutely least area (or least weight) essential torus must have the least number of double curves among all the essential tori.  In particular, an absolutely least area (or least weight) torus has only one or two double curves and has the 4-plane and 1-line properties as described in \cite{S}.
\end{theorem}

Since there are only finitely many normal surfaces (up to normal homotopy) with weight no more than $W$, we can enumerate all the normal tori with weights bounded by $W$.  By Theorem~\ref{T:LN}, if $M$ is a Seifert fiber space of hyperbolic type, one of these normal tori has the 4-plane and 1-line properties.  If $T$ is a torus with the 4-plane and 1-line properties, then as in \cite{S, HS}, one can eliminate all triple points by performing some simple homotopies.  In fact, for any essential torus $g:T\to M$ with the 4-plane and 1-line properties as described in \cite{S}, if there are triple points, then by the argument in Lemma 4.4 of \cite{S}, the closure of at least one component of $M-g(T)$ must be an embedded 3-side ``football region", i.e., an embedded 3-ball $B^3$ whose boundary consists of 3 bigons and $int(B^3)\cap g(T)=\emptyset$ (a bigon is a disk in $g(T)$ whose two vertices are triple points and whose two edges are in the double curves of $g(T)$).  By moving one bigon across the other two bigons in a small neighborhood of the ``football region", we can cancel a pair of triple points, and eventually such homotopy can eliminate all triple points.  

So, we perform these homotopies on each of the normal tori with weight less than $W$.  If we cannot find such a ``football region" and get a torus without triple points by such homotopies, $M$ is not a Seifert fiber space of hyperbolic type.  If we find a torus without triple points, by checking whether the complement of this immersed torus is a union of solid tori, we can determine whether $M$ is a Seifert fiber space.

\begin{remark}
\begin{enumerate}
\item The key reason for the existence of such a nice ``football region" is the 4-plane and 1-line properties.  For immersed least area surface in general, there may not exist such an embedded ``football region".  Also, if there are pieces of the immersed surface in the interior of a ``football region", such a homotopy may increase the number of the triple points.  In fact, there are examples of immersed surface that one cannot reduce the number of the triple points without increasing the number of the triple points first.
\item  Theorem~\ref{T:LN} is not easy to prove.  Although Theorem~\ref{T:LN} makes the algorithm faster, this theorem is not absolutely essential for our algorithm.  The following much weaker version of Theorem~\ref{T:LN} is fairly easy to prove (by cutting and pasting on multiple copies of $f(T)$), and is sufficient to make the algorithm work.  Given an essential least weight immersed torus $f:T\to M$, where $M$ is a non-Haken small Seifert fiber space of hyperbolic type, the weight of Scott's tori (those with one or two double curves and the 4-plane and 1-line properties) are less than $h(w, d)$, where $h$ is an explicit function of $w=weight(f(T))$ and $d$ that is the number of the double curves of $f(T)$.
\end{enumerate}
\end{remark}

\end{psfrags}

\end{document}